\documentclass[aps,pra,notitlepage]{revtex4-1}
\usepackage[utf8]{inputenc}
\usepackage{amsmath}
\usepackage{tikz}
\newcommand{\ie}{,~i.e.,}
\begin{document}
\title{Monty Hall problem revisited once more}
\author{Francisco A. B. Coutinho}
\author{Eduardo Massad}
\affiliation{Faculdade de Medicina, University of S\~ao Paulo}
\author{Luiz N. Oliveira}
\affiliation{Instituto de Física de S\~ao Carlos, University of S\~ao Paulo, \\
  C.~P.~369 --- 13560 São Carlos, SP, Brazil}
\begin{abstract}
  The Monty Hal problem is an attractive puzzle. It combines simple
  statement with answers that seem surprising to most audiences. The
  problem was thoroughly solved over two decades ago. Yet, more recent
  discussions indicate that the solution is incompletely
  understood. Here, we review the solution and discuss pitfalls and
  other aspects that make the problem interesting.
\end{abstract}
\maketitle

The Monty Hall problem (MHP), equivalent to the three-prisoner puzzle,
is a brain teaser from probability theory. In 1990, vos Savant
succinctly presented a solution \cite{1990Sav9}, which was spelled out
by the same author in two subsequent notes \cite{1990Sav2,1991Sav17}.
Shortly afterwards, her solution was challenged by Gillman
\cite{1991Gil8,1992Gil3}, who argued that she had implicitly
considered a modified version of the game, and proceeded to presenting
a thorough solution. Nonetheless, the lucid explanations in Gillman's
articles seem not to have been thoroughly absorbed, since more recent
discussions have presented incomplete analyses of the problem. Very
recent examples, which we single out because they appear in books of
outstanding quality, are due to Ben Naim \cite{2015Nai} and Miller
\cite{2017Mil}. Yet another survey of the various aspects of MHP seems
therefore warranted.

\section{The Monty Hall game}
\label{sec:monty-hall-game}

\begin{figure}[!ht]
  \centering
  \begin{tikzpicture}\label{fig:1}
    \draw(-2.75,-2.75)rectangle(2.75,2.75);
    \draw(0,0)circle(2cm);
    \foreach \x/\y in {0/2, -1.73/-1, 1.73/-1}{
      \begin{scope}[xshift=\x cm, yshift= \y cm]
        \fill[green!30!white](-0.5,-0.5)rectangle(0.5,0.5);
      \end{scope}
      }

    \draw[red](0,2)node{$T$}; 
    \draw[red](-1.73,-1)node{$L$}; 
    \draw[red](1.73,-1)node{$R$}; 
  \end{tikzpicture}
  \caption[Game]{Positions of the three boxes after the contender has stated
    her initial choice. The structure has been rotated so that the
    chosen box lies at the top (T). The other two boxes are now at the
    positions labeled L and R.}
\end{figure}
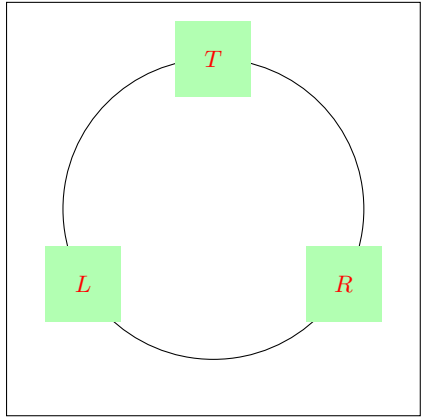

The rules of the game are simple. Monty Hall shows three closed boxes
to a contender\textemdash call her Portia. One of them contains an
automobile, while a goat stands inside each of the other two. The
boxes are symmetrically disposed around a circular track, which can be
freely rotated around its center. Portia chooses one of the boxes,
and the track is immediately rotated to the disposition in
Fig.~\ref{fig:1}.

Monty Hall knows where the car is. He opens one of the two boxes at the bottom of the
figure\textemdash box $L$, for definiteness\textemdash to show Portia that
it contains a goat, and asks her whether she wants to stick to her
choice or switch to box $R$. The problem asks for the probability
$P_R$ that Portia will win the car if she decides to switch. 

\section{Gillman's versions}
\label{sec:gillmans-version}
Gillman pointed out that the game admits two variants
\cite{1991Gil8,1992Gil3}. One of them, which he calls version~I,
follows the rules in Section~\ref{sec:monty-hall-game}.

\subsection{Gillman's Game I}
\label{sec:gillmans-game-i}

Gillman observes that a complete statement of the game in
Section~\ref{sec:monty-hall-game} calls for knowledge of Monty Hall's
tactics. The host has some freedom. If the car lies in box $L$ ($R$),
in Fig.~\ref{fig:1} then his only option is to open door $R$ ($L$). If
the car lies in box $T$, however, the host may open either door. We
will call $q$ be the probability that he will open door $L$. To
compute her odds, Portia must know $q$. 

Let us assume that she does know $q$ and consider two illustrative
situations: $q=0$ and $q=1$.  With $q=0$, Monty Hall will blindly
open door $R$ whenever given the choice. By opening door $L$, he in
practice admits to Portia that he had no choice, which means that the
car must be behind door $R$. Portia, therefore, asks to switch with the
winner's smile on her lips. In this case, the probability of winning
by switching is $P_R=1$.

With $q=1$, by contrast, if the car is in box $T$ Monty Hall will open
box $L$. If the car is in box $R$, he will also open box $L$. Hence,
the only information conveyed to Portia when he opens $L$ is the
obvious one: that box $L$ holds a goat. To switch or not to switch is
then Portia's question, because the chances are equal. In other words,
$P_R=1/2$ for $q=1$. 

We can see that the probability depends on $q$. To come to the same
conclusion from another perspective, consider the probability
$P_T$, of winning without switching. Clearly, $P_T=1-P_R$. Before
Monty opened the door, the probability was $P_T=1/3$. When he opens
door $L$, a possibility is eliminated, but this will not change $P_T$
if boxes $L$ and $R$ are equivalent, i.e., if Monty is unbiased
towards either box, i.e, if $q=1/2$. We conclude that, with $q=1/2$,
the probability of winning by not switching is $P_T=1/3$, and the
probability of winning by switching is $P_R=2/3$.


With $q=1/2$, therefore, Savant's argument applies, and the
probability of winning by switching is $P_R=2/3$. Clearly, vos Savant
had the unbiased game in mind when she wrote on the subject
\cite{1990Sav9,1990Sav2,1991Sav17}, i.e., she considered the $q=1/2$
special case of Gillman's Game I. Gillman argued that she had another
game in mind. That interpretation seems unfair, even though her
solution does apply to Gillman's Game II, which we describe next.

\subsection{Gillman's Game II}
\label{sec:gillmans-game-ii}

In the second version of the game, Portia must announce her decision to switch or not to switch before Monty opens the door. Knowing the host's bias is now of no help to Portia, since she must make her decision before Monty opens either box. In that case, von Savant's argument is impeccable, for Portia's choice cannot affect the probability $P_T=1/3$ that the car be in the top box. The probability that Portia wins by switching is therefore $2/3$.
\section{Formal analysis}
\label{sec:formal-analysis}
We now present mathematical proofs of the results derived on
intituitive grounds in Secs.~\ref{sec:gillmans-version}~and \ref{sec:gillmans-game-i}.
\subsection{Game I}
\label{sec:game-i}

\subsubsection{First solution}
We ask for the conditional probability that the car lie in box $i$
($i=L,R$), given that Monty has opened box $j$ ($j=L,R$). In other
words, we ask for the conditional probability $P(C_j|H_i)$, where
$C_j$ denotes \emph{event $j$}\ie\ the car is in box $j$, and $H_i$
denotes event $i$\ie\ Monty Hall has opened box $i$. As Portia
knows, if the car is in box $T$, the probability that Monty will open
box $R$ is
\begin{align}
  \label{eq:4}
  P(H_R|C_T) = q.
\end{align}

Let us now turn to Bayes's formula, which states that
\begin{align}
  \label{eq:1}
  P(H_R)P(C_L|H_R)= P(C_L)P(H_R|C_L),
\end{align}
and, likewise, that
\begin{align}
  \label{eq:2}
  P(H_R)P(C_T|H_R) = P(C_T)P(H_R|C_T).
\end{align}

The \emph{a priori} probabilities that the car be in boxes $L$ or $T$
are identical: $P(C_L)=P(C_T)=1/3$.
Division of Eq.~\eqref{eq:2} by Eq.~\eqref{eq:1} therefore shows that
\begin{align}
  \label{eq:3}
  \dfrac{P(C_L|H_R)}{P(C_T|H_R)} = \dfrac{P(H_R|C_L)}{P(H_R|C_T)}. 
\end{align}

Moreover, $P(H_R|C_L)=1$, because Monty Hall is forced to open box $R$
when the car is in box $L$. From Eqs.~\eqref{eq:4}~and \eqref{eq:3} it follows
that
\begin{align}
  \label{eq:5}
  \dfrac{P(C_L|H_R)}{P(C_T|H_R)} = \dfrac1q.
\end{align}

We then recall that events $C_T$ and $C_L$ are mutually
exclusive when event $H_R$ takes place (in other words, when Monty
opens door $R$, he is telling Portia that the car is either in box $T$
or in box $L$). This shows that the numerator and denominator of the
fraction on the left-hand side of Eq.~\eqref{eq:5} add up to unity:
\begin{align}
  \label{eq:6}
  P(C_L|H_R)+P(C_T|H_R) = 1.
\end{align}

Solution of the system of Eqs.~\eqref{eq:5}~and \eqref{eq:6} yields
the desired result:
\begin{align}
  \label{eq:7}
  P(C_L|H_R) = \dfrac1{1+q}.
\end{align}

As expected from our discussion in Sec.~\ref{sec:gillmans-game-i}, the
probability that Portia win the car by switching therefore varies 
from unity, for $q=0$, to $1/2$, for $q=1$. When Monty is unbiased,
$q=1/2$, and the probability is $2/3$, as von Savant predicted.

\subsubsection{Second solution}
\label{sec:second-solution}
We now follow Isaac \cite{1995Isa} to present an alternative solution
that seems more attractive because it starts by surveying the space
$\mathcal{S}$ of possible events. Let $C_i$ ($i=T,L,R$) denote the
event in which the car lies in box $i$. Let $S_k$ ($k=L,R$) denote the event that
Portia switches to box $k$, let $W$ indicate that Portia wins\textemdash
she drives the car home\textemdash and let $L$ indicate that Portia
loses\textemdash she trails back home in the company of a goat.
Thus, for instance, the event $C_R$ is incompatible with $S_L$: if the
car is behind door $R$, then Monty is forced to open door $L$, and
Portia would not switch to the open box, inside which she can see a goat.

The sample space comprises four combinations of events:
\begin{align}
  \label{eq:8}
  \mathcal{S} = \{(C_T,S_R,L), (C_T, S_L, L), (C_L, S_L,W), (C_R,S_R,W)\}, 
\end{align}

We can now compute the probability for each combination on the
right-hand side of Eq.~\eqref{eq:8}. Each individual event $C_i$
($i=T,L,R$) occurs with probability $1/3$. The combined event
$(C_T,S_R,L)$ can only occur when, given that the car is in box $T$,
the host chooses to open box $L$. Its probability is, therefore,
$(1/3)(1- q)$. Likewise, the combination $(C_T,S_L,L)$ will only
occur if the host chooses to open box $R$. Its probability is
$(1/3)q$. Finally, the last two combinations, which are
independent of Monty's bias, have probability $1/3$ each.

When the host opens door $R$ (event $H_R$), Portia rules out the first and
the last combinations of events on the right-hand side of
Eq.~\eqref{eq:8}. She can then compute her chance of winning by
switching from the probabilities for the other two combinations:
\begin{align}
  \label{eq:9}
  P_{S_L|H_R} = \dfrac{\dfrac{1}3}{\dfrac{q}3+\dfrac{1}{3}}, 
\end{align}
which yields the desired result
\begin{align}
  \label{eq:10}
  P_{S_L|H_R} = \dfrac{1}{1+q},
\end{align}
in agreement with Eq.~\eqref{eq:7}

\subsection{Game II}
\label{sec:game-2}
Gillian's Game II is somewhat different. As in Game I, Portia chooses a box, which is promptly moved to the top position, as in Fig.~\ref{fig:1}. Monty then announces that he will soon show her a goat inside one of the two boxes at the bottom in Fig.~\ref{fig:1}. Immediately after the announcement, before opening any door, he asks her whether she will switch to the other box. If Portia chooses to switch, what are her chances of winning the car?

\subsubsection{First solution}
\label{sec:first-solution}

As before, let us call $C_i$ ($i=T,L,R$) the event that the car is in
position $i$. As in Sec.~\ref{sec:second-solution}, we can see that
the probability $P(C_i)=1/3$ ($i=T,L,R$). Call $P(W|C_i)$ the
probability that Portia wins by switching if the car is in box $i$
($i=T, L,R$). The probability of winning by switching is
\begin{align}
  \label{eq:11}
P(W)= P(W|C_T)P(C_T) + P(W|C_L)P(C_L)+P(W|C_R)P(C_R).
\end{align}

If the car is inside box $T$, switching will make Portia lose. Hence
$P(W|C_T)=0$. She knows that, in this case, the probability that Monty will choose box $R$ (box $L$)
is $q$ ($1-q$), but that information is of no help, because the probability of
winning is zero, either way.

On the other hand, $P(W|C_L)=1$, for if the car is in box
$L$, Monty must open box $R$, and Portia will be compelled to make the right
choice\ie\ open box $L$. Likewise, $P(W|C_R)=1$. It follows from
Eq.~\eqref{eq:11} that
\begin{align}
  \label{eq:12}
P(W)=   0\times\dfrac13+1\times\dfrac13+1\times\dfrac13,
\end{align}
that is,
\begin{align}
  \label{eq:13}
  P(W) = \dfrac{2}3.
\end{align}

By contrast, if Portia decides to stick with her initial choice,
$P(W|C_T)$ will be unitary, while $P(W|C_L)=P(W|C_R)=0$. In this case,
Eq.~\eqref{eq:12} yields $P(W)=1/3$.

\subsubsection{Second solution}
\label{sec:second-solution-1}

As in Sec.~\ref{sec:second-solution}, let us follow Isaac \cite{1995Isa} and consider the four possible combinations of events. The sample space is again given by Eq.~\eqref{eq:8}, and the probability that each combination occurs is the one we have already computed: $(1-q)/3$ for $(C_T, S_R,L)$, $q/3$ for $(C_T, S_L,L)$, $1/3$ for $(C_L, S_L,W)$, and $1/3$ for $(C_R, S_R,W)$. In either of the latter two events, Portia drives the car home. The probability of winning by switching is therefore the sum of the two probabilities, and we recover Eq.~\eqref{eq:13}.

Equation~\eqref{eq:13} agrees with the result derived in Refs.~\onlinecite{1991Sav17,1990Sav9,1990Sav2}. This does not mean, however, that von Savant was playing Game II. As already pointed out, all evidence indicates that she was considering Game 1 with $q=1/2$\ie\ with an unbiased host. Under that condition Eqs.~\eqref{eq:7}~and \eqref{eq:10} are also equivalent to Eq.~\eqref{eq:13}.

\section{A long series of Game I}
\label{sec:long-series-game}

While discussing Game~I, In Sec.~\ref{sec:game-i}, we assumed that Monty had opened door $R$, and our computation of the probability relied on that assumption. If the game is played not once, but numerous times, we cannot assume that box $R$ will be opened every time. In a long series of games, the host's bias becomes irrelevant. The distinction between Games I and II is washed out, as the following discussion shows.

In a long sequence, the car will be in box $T$ in $1/3$ of the events. In those events, Monty has his way. Whether he has or has not a preference for one of the doors is entirely imaterial, however, for in the end Portia will take a goat home, if she decides to switch. Of course, if he is biased towards one of the doors, Portia will more often take home the goat in other box, but that will be of no comfort to her.

In the other $2/3$ of the events, the car will lie either in box $R$ or in box $L$. In those events, Portia's decision to switch will add a brand new automobile to her assets. The odds of winning by switching, therefore, are $2/3$, just as if Monty and Portia were playing Game II.

\bibliography{mHall}
\end{document}